\documentstyle{article}
\newtheorem{df}{ \sc Definition}[section]

\newtheorem{ex}[df]{ \it Example}
\newtheorem{pr}[df]{ \sc Proposition}

\newtheorem{th}[df]{ \sc Theorem}
\newtheorem{cor}[df]{ \sc Corollary}
\newtheorem{re}[df]{ \it Remark}

\newtheorem{lem}[df]{ \sc Lemma}

\newtheorem{eq}[df]{\rm}
\def\barr{\overline}

\def\p{\phantom{\mpd{}^\mpd{}_\mpd{}}}
\def\mpr#1{\;\smash{\mathop{\hbox to 20pt{\rightarrowfill}}\limits^{#1}}\;}
\def\epi#1{\;\smash{\mathop{\hbox to 20pt{\rightarrowfill}\hskip
-13pt\rightarrow}\limits^{#1\,}}\;}
\def\epii{\smash{\mathop{\hbox to 14pt{\rightarrowfill}\hskip
-10pt\rightarrow}}}
\def\mono{\lhook\joinrel\relbar\joinrel\rightarrow}
\def\mpL#1{\;\smash{\mathop{\hbox to 20pt{\leftarrowfill}}\limits_{#1}}\;}
\def\mpl#1{\;\smash{\mathop{\hbox to 20pt{\leftarrowfill}}\limits^{#1}}\;}
\def\mpd#1{\big\downarrow\rlap{$\vcenter{\hbox{$\scriptstyle#1$}}$}}
\def\mdp#1{\llap{$\vcenter{\hbox{$\scriptstyle#1$}}$}\big\downarrow}

\font\Ch=msbm10
\def\Rt{\hbox{\Ch R}}
\def\Zt{\hbox{\Ch Z}}

\def\Qt{\hbox{\Ch Q}}
\def\Ct{\hbox{\Ch C}}
\def\Pt{\hbox{\Ch P}}
\def\Ft{\hbox{\Ch F}}

\def\Proof{\noindent{\it Proof. }}
\def\qed{\hfill$\Box$\vskip10pt}

\begin{document}

\author{Andrzej Weber\footnote{Supported by KBN 2P03A 00218 grant.
I thank Institute of
Mathematics, Polish Academy of Science for hospitality.} }
\date{}
\title{Weights in the cohomology of toric varieties}
\maketitle

\begin{abstract}

\vskip5pt
We describe the weight filtration in the cohomology of toric
varieties. We present a role of the Frobenius
automorphism in an ele\-mentary way.
We prove that equivariant intersection homology of an arbitrary toric
variety is pure.
We obtain results concerning Koszul duality: nonequivariant
intersection cohomology is equal to the cohomology of the Koszul
complex $IH_T^*(X)\otimes H^*(T)$.
We also describe the weight filtration in $IH^*(X)$.
\end{abstract}

\section{Introduction} Let $X$ be a smooth toric variety. If $X$ is
complete then its cohomology coincides with the Chow ring $A^*(X)$.
Therefore the Hodge structure is not very interesting:
$H^{k,k}(X)\simeq A^k(X)\otimes\Ct$ and $H^{k,l}(X)=0$ for $k\neq l$.
If $X$ is not complete, then the cohomology of $X$ is equipped with
the weight filtration constructed by Deligne, \cite{D2}:
$$W_{0}H^k(X)\subset\dots \subset
W_{2k}H^k(X)=H^k(X)\,.$$
Since $X$ is smooth $W_{k-1}H^k(X)=0$. In the toric case
the weight filtration has the property:
$$W_{2l}H^k(X)=W_{2l+1}H^k(X)$$ for each $l$. The pure
Hodge structure on
$$Gr_{2l}^WH^k(X)=W_{2l}H^k(X)/W_{2l-1}H^k(X)$$
consists only of the Hodge type $(l,l)$.

If $X$ is singular we replace ordinary cohomology by intersection
cohomology. Since $X$ can be given a locally conical structure in a
metric sense, the intersection cohomology may be interpreted as
$L^2$-cohomology of the nonsingular part, see \cite{Ch, CGM}.
Nevertheless we avoid to talk about the $L^2$--Hodge structure. We are
interested in the weight filtration, which is defined via
reduction of $X$ to a finite characteristic field $\Ft_p$, see
\cite{BBD}. The weight filtration is the
filtration coming from the eigenvalues of the Frobenius
automorphism acting on the \'etale intersection cohomology.

Fortunately, for toric varieties the Frobenius automorphism of
intersection cohomology can be induced from an endomorphism of the
complex points of $X$. This map preserves orbits and it is covered
by a map of intersection cohomology sheaves. Finally, the group
$IH^k(X)$ is decomposed into the direct sum of the eigenspaces
with eigenvalues $p^l$, $l=\lceil {k\over 2} \rceil ,\lceil
{k\over 2} \rceil+1,\dots, k$. We stress that although the concept of
weights is highly nontrivial and abstract, for toric varieties
whole theory reduces to easy computations involving an action of
a down-to-earth map.

On the other hand we consider
equivariant cohomology taken with respect to the big torus $T=(\Ct^*)^n$
acting on $X$, \cite{BBFK1}. If $X$ is singular we prefer to replace
usual cohomology by the intersection cohomology.
We prove that the equivariant intersection cohomology of $X$ is pure
even if $X$ is not complete. This means that the only eigenvalue which can
occur in $IH^{2k}_T(X)$ is $p^k$, whereas $IH^{2k+1}_T(X)=0$.

Koszul duality allows to recover usual intersection
cohomology from equivariant one \cite{GKM, MW, Fr1}. In general one has to 
know
not
only the cohomology groups, but also whole complex in the derived
category of $H^*(BT)$-modules. One has a spectral sequence with
$E^{k,l}_2=IH^k_T(X)\otimes H^l(T)$ converging to $IH^{k+l}(X)$.
The differential in the $E_2$ table is the usual Koszul differential:
We identify $H^*(T)$ with the exterior algebra of the dual of the Lie
algebra $\bf t^*$, and the cohomology of $BT$ with the symmetric
algebra. The differential $d_{(2)}=d_{Koszul}$ has the form:
$$d_{(2)}(x\otimes \xi)=\sum_{j=1}^n x\lambda^j\otimes
i_{\lambda_j}\xi\,, \quad{\rm for}\; x\in IH^*_T(X)\,,\;\;\xi\in\Lambda{\bf
t}^*\,.$$
Here $\{\lambda_j\}$ is a basis of $\bf t$,
the elements $\lambda^j$ of the dual basis are generators of
$H^*(BT)=S{\bf  t}^*$,  and  $i_{\lambda_j}$   stands   for   the 
contraction.
We prove using weight argument that all the higher differentials of the
spectral sequence vanish if $X$ is toric.
Therefore the nonequivariant cohomology of $X$ are the cohomology of
the Koszul complex.
Moreover, the weight filtration coincides with the filtration given
by the spectral sequence. The Koszul complex itself splits into a
direct sum of subcomplexes, each computing the graded piece of the
weight filtration $Gr^W_*IH^*(X)$.
We note that the nonequivariant intersection cohomology of $X$ is pure if and
only
if the equivariant intersection cohomology is free over $H^*(BT)$.
The properties of the weight filtration are reflected by the
behavior of Poincar\'e polynomials, which are the weighted Euler 
characteristics.

Cohomology and intersection cohomology of toric varieties have
attracted surprisingly many authors. Jurkiewicz (\cite{Ju})
and Danilov (\cite{Da})
have computed cohomology of complete smooth
toric varieties. Next we would have to present a long list of
papers. Instead we suggest to check the references e.g. in \cite{BP}.
A complex computing intersection cohomology of a toric variety
was described in \cite{Is}.
Our paper is based on the approach of \cite{Fi}, [BBFK1-2]
combined with \cite{DL}. We will assume that the reader is
familiar with basic theory of toric varieties (\cite{Fu}),
intersection cohomology (\cite{GM}) and equivariant intersection
cohomology (\cite{Br}). Our goal is to expose the role of the
weight filtration. Toric varieties serve as a toy models. The
reader who is not familiar with intersection cohomology may
replace it by usual cohomology and assume that $X$ is defined by a
simplicial fan, i.e.~$X$ is rational homology manifold.

Now, we would like to explain terminology
concerning formality which plays an important role in our consideration.
\begin{enumerate}
\item A manifold is called formal if the algebra of differential
forms $\Omega^*(X)$ is quasiisomorphic to its cohomology {\it as
a dg-algebra.} In the definition of formality for an arbitrary
topological space the algebra of forms is replaced by the
Sullivan--de Rham
complex. 
Except from K\"ahler manifolds the classifying space of a
connected Lie group $BG$ is an example of a formal space. 

\item Suppose $B$ is a formal space and $X\to B$ is a map. Then
$\Omega^*(X)$ is a module over $\Omega^*(B)$. A natural notion of
formality over $B$ would be the demand that $\Omega^*(X)$ is
quasiisomorphic to its cohomology {\it as a dg-module over
$\Omega^*(B)\simeq H^*(B)$.}  Formality of $ET\times_TX$ over
$BT$ implies that 
$$\matrix{ H^*(X)&=&H^*(\Omega^*(ET\times_TX)\otimes
\Lambda^*,d_{Koszul})\hfill\cr 
&=& H^*(H^*(ET\times_TX)\otimes \Lambda^*,d_{Koszul})\,.\cr}$$
This is exactly the content of our Theorem \ref{kosz} for simplicial
toric varieties.
Recently M. Franz \cite{Fr2} have shown that $ET\times_TX$ is formal
over $BT$ (even with $\Zt$ coefficients) for smooth toric varieties.

\item If $B=BG$ and $\Omega^*(EG\times_GX)$ is quasiisomorphic to its
cohomology {\it as an algebra over $H^*(BG)$} then $X$ is called
$G$-formal in \cite{Li}.  Smooth toric varieties
are formal in the above sense by \cite{NR}.

\item The notion of {\it equivariantly formal} space was introduced
in \cite{GKM}. Before it was called
{\it totally nonhomologous to zero}. The name does not fit to the 
scheme of the previous
definitions. It just means that $\Omega^*(EG\times_GX)$ is free (up
to a quasiisomorphism)
dg-module over $H^*(BG)$. It is equivalent to the statement that
$H^*_G(X)$ is a free module over $H^*(BG)$.
\end{enumerate}

I would like to thank Matthias Franz for valuable remarks and
comments.

{\def\contentsname{\normalsize\bf Contents:}\tableofcontents}

\section{Mixed Hodge structure}

Let $X$ be a smooth possibly noncomplete algebraic variety. According
to Deligne (\cite{D2}) one defines an additional structure on the
cohomology of $X$. One finds a completion $X\subset \barr X$, such
that $\barr X\setminus X=\sum_{i=1}^\alpha D_i$ is a smooth divisor with
normal crossings. For $0\leq k\leq\alpha$ let $D^{(k)}$ denote
the disjoint union of the $k$-fold intersections of the
components of the divisor $D$ and set $D^{(0)}=\barr X$. Deligne has
constructed a spectral sequence with
\begin{eq}\hfil$E^{k,l}_1=H^{2k+l}(D^{(-k)})\Rightarrow
H^{k+l}(X)\label{ss}$\end{eq}
(the coefficient are in $\Ct$).
This spectral sequence degenerates on $E_2$. The limit filtration is
the weight filtration. The quotients are equipped with the pure Hodge
structure, they are decomposed into the $(p,q)$ summands:
$$Gr^W_kH^l(X)=\bigoplus_{p+q=k} {\cal H}_{(l)}^{p,q}\,.$$
The groups $Gr^W_kH^l(X)$ are quotients of subgroups of
$H^{l-2k}(D^{(k)})$.

If $X$ is a smooth toric variety then there exists a smooth toric
variety compactifying $X$. It can be chosen so the divisors at
infinity are smooth toric varieties as well as each component of
$D^{(k)}$. The resulting weights can only be even. Moreover:

\begin{pr} For a smooth toric variety each $Gr^W_{2k}H^l(X)$ is of the
type $(k,k)$ and $Gr^W_{2k+1}H^l(X)=0$.
\end{pr}

\begin{re}\rm There is an easy method of constructing a spectral sequence
converging to the weight filtration. Just consider the Leray spectral
sequence of the inclusion $X\subset \barr X$.
Then $E^{k,l}_2=H^k(D^{(l)})$. The limit filtration has to be
shifted in order to have $W_lH^l(X)=im(H^l(\barr X)\rightarrow
H^l(X))$. This spectral sequence degenerates on $E_3$. In fact, up to
a renumbering it is isomorphic to the Deligne spectral sequence.\end{re}

The construction of Deligne was motivated by the previous work on
Weil conjectures, \cite{D1}. One considers varieties defined over a finite
field. Instead of the Hodge structure one has an action of the
Frobenius automorphism. The eigenvalues on $H^l_{et}(X)$
can have absolute values equal to $p^{k \over 2}$ with $k=l,l+1,\dots, 2l$
if $X$ is smooth.
Each complex variety is in fact defined over a finitely
generated ring; toric varieties are defined over $\Zt$.
The \'etale cohomology of the variety reduced to the finite
base field is isomorphic to the usual cohomology for almost
all reductions. In the toric case every reduction $\Zt\rightarrow \Ft_p$
is good.
The weight filtration coincides with the filtration
of the \'etale cohomology:
\begin{eq}\hfil\label{dfil}$W_kH^l_{et}(X)=\bigoplus_{|\lambda|\leq
p^{k\over 2}}V_\lambda\,,$\end{eq}
where $V_\lambda\subset H^l_{et}(X)$ is the eigenspace of the
eigenvalue $\lambda$. The eigenvalues
 appearing here belong to $\overline
{\Qt}_\ell\simeq \Ct$ in general, but in the toric case $\lambda$ is a
power of $p$.

\section{Frobenius automorphism of toric varieties}

As noticed by Totaro (\cite{To}), the toric varieties are rare examples of
complex varieties admitting an endomorphism which coincides with
the Frobenius automorphism after reduction.

Let $p>1$ be a natural number. Raising to the $p$-th power is a map
of the torus
$\psi_p:T\rightarrow T$. There exists an extension
$\phi_p:X\rightarrow X$, such that the following diagram commutes:
\begin{eq}\label{fipsi}\hfil$\matrix{T\times X&\mpr{\mu}& X\cr
\mdp{\psi_p\times\phi_p}& &\mpd{\phi_p}\cr
T\times X&\mpr{\mu}& X&.\cr}$\end{eq}
Here $\mu$ is the action of $T$.
The map $\phi_p$ is given by a subdivision of the lattice. 
It has been recently applied to study equivariant
Todd class in homology of singular toric variety in \cite{BZ}

The following observations (already made in \cite{To})
illustrate the nature of the weight
filtration very well. To prove them one applies elementary properties
of cohomology. 

\subsection{$X$ complete and smooth}

Let $X$ be smooth and complete toric variety. Then every homology
class is represented by the closure of an orbit. The restriction of
$\phi_p$ to a $k$-dimensional orbit is a cover of the degree $p^k$.
Therefore the induced action of $\phi_p$ on $H_{2k}(X)$ is the
multiplication by $p^k$. The conjugate action on $H^{2k}(X)$ is
the multiplication by $p^k$. We say then that the cohomology of $X$
is pure (see Definition \ref{pude} below).

\subsection{$X$ smooth}

If $X$ is smooth, but not necessarily complete, then one has Deligne
spectral sequence \ref{ss}. The induced action of $\phi_p$ on
$Gr^W_{2k}H^l$ is the multiplication by $p^k$. It can be verified by
comparison with the \'etale cohomology or in a elementary way: if $Y$ is
the closure an orbit of the codimension one, then
the following diagram commutes
$$\matrix{H^*(Y)& \mpr{i_!} & H^{*+2}(\barr X)\cr
\mdp{\phi_p^*}&\p&\mpd{\phi_p^*}\cr
H^*(Y)& \mpr{pi_!} & H^{*+2}(\barr X)&.\cr}$$
Here $i_!$ is Gysin the map induced by the inclusion via Poincar\'e
Duality. (It is a general rule that Gysin map shifts the weight by
the codimension of the subvariety.)
The differential in the Deligne spectral sequence \ref{ss} are induced by the
inclusions. Therefore, to have an equivariant spectral sequence with
respect to $\phi_p$ one encodes the action in the usual way:
$$E^{-k,l}_1=H^{2k+l}(D^{(k)})(k)\,,$$
where the symbol $(k)$ denotes action of $\phi_p$ via the
multiplication by $p^k$.
The spectral sequence degenerates on $E_2$. This is because the
higher differentials
$$d_{(i)}:E^{k,l}_i\rightarrow E^{k+i,l-i+1}_i$$
do not preserve the eigenspaces of $\phi_p^*$. The domain of $d_{(i)}$
has the eigenvalue $p^{l\over 2}$ and the target has $p^{l-i+1\over
2}$ (the relevant values of $l$ and $l-i+1$ are even).
For the limit filtration of $H^*(X)$ we have:

\begin{pr} The action of $\phi_p^*$ on $Gr^W_{2k}H^l(X)$ is the
multiplication by $p^k$.\end{pr}

Again this is just a simple characterization 
of the weight filtration for toric varieties. On the other hand,
not appealing to the general theory, we can define weights in the
following way.
To be consistent with the usual terminology we will say that:

\begin{df} A vector space $V$ is of {\it weight} $k$ if
\begin{itemize}
\item  $\phi_p$ acts on $V$ as
multiplication by $p^{k/2}$ for $k$ even,
\item  $V=0$ for $k$ odd.
\end{itemize}
\end{df}

\begin{df}\label{pude}
The cohomology group $H^k(X)$ is {\it pure} if it is of weight $k$.
\end{df}

It turns out (see Corollary \ref{purpur}) that our definition of
purity agrees with the one in \cite{BBFK2}, \S2. Namely the
condition on the eigenvalues in even cohomology redundant. It is implied by
vanishing of odd cohomology.

In fact the decomposition of $H^*(X)$ into the eigenspaces of
$\phi_p^*$ splits the weight filtration into a gradation, \cite{To}.

\begin{ex}\label{prz} \rm Let
$X=\Pt^1\times\Pt^1\setminus\{(0,\infty),(\infty,0)\}$.
Let $\barr X$ be the compactification of $X$ obtained by blowing up
two removed points in $\Pt^1\times\Pt^1$.
The Deligne spectral sequence has the following $E_1$ table
$$\begin{tabular}{ccc|c}
0&$\Ct^2(2)$&$\Ct(2)$&4\\
0&0&0&3\\
0&$\Ct^2(1)$&$\Ct^4(1)$&2\\
0&0&0&1\\
0&0&$\Ct$&0\\ \hline
-2&-1&0&\\
\end{tabular}\qquad d_{(1)}:E^{k,l}_1\rightarrow E^{k+1,l}_1.$$

\noindent The second table is equal to $E_\infty$:
$$\begin{tabular}{ccc|c}
0&$\Ct(2)$&0&4\\
0&0&0&3\\
0&0&$\Ct^2(1)$&2\\
0&0&0&1\\
0&0&$\Ct$&0\\ \hline
-2&-1&0&\\
\end{tabular}
\qquad \phantom{d_{(1)}:E^{k,l}_1\rightarrow
E^{k+1,l}_1.}$$
The cohomology of $X$ is following:
$$ \Ct,\,0,\,\Ct^2(1),\,\Ct(2),\,0$$
in the dimensions 0,1,2,3,4. The group $H^3(X)$ is not pure, since it is
of the weight 4.
\end{ex}

\subsection{$X$ arbitrary}
For singular $X$ we study intersection cohomology (\cite{GM}) instead
of usual cohomology. The map $\phi_p$ preserves orbits, which form a
stratification of $X$. Therefore $\phi_p$ induces a map of
intersection cohomology groups $IH^*(X)$.
Moreover this map exists on the level of sheaves. Thus it is natural
with respect to inclusions. Avoiding the formalism of \cite{BBD}
we can define a weight filtration by the formula \ref{dfil}. Later it
will be clear, that the
action of $\phi_p^*$ decomposes $IH^*(X)$ into a sum of eigenspaces
with the eigenvalues $p^k$ as in the smooth case, see \ref{kosz}.

\section{Equivariant cohomology}

Although the equivariant cohomology and equivariant intersection
cohomology are complicated objects, they have surprisingly nice properties
for toric varieties. Equivariant cohomology is defined by means of
Borel construction:
$$H^*_T(X)=H^*(ET\times_TX)\,.\qquad
IH^*_T(X)=H^*(ET\times_TX;IC^{_{^\bullet}}_T)\,.$$
For the basic features we refer
to \cite{Br} and \cite{BBFK1}.

\subsection{Smooth case} \label{poly}

Equivariant cohomology of a smooth (or
just simplicial) toric variety can be easily recovered from the fan.
If $\Delta$ is a simplicial fan, then the associated
toric variety $X_\Delta$ is a rational homology manifold. The
equivariant cohomology can be identified with the algebra of the
continuous functions on the support of the fan, which are polynomial
on each cone of $\Delta$. The odd part of the equivariant cohomology
vanishes also for noncomplete fans. From the commutativity of \ref{fipsi}
it follows that $\phi_p$ acts on the equivariant cohomology.

We will prove a key theorem:

\begin{th} \label{pec} The equivariant cohomology of a smooth toric
variety is pure.\end{th}

\begin{re}\rm Since the space $ET$ is a limit of algebraic varieties
(as defined in \cite{Br})
or a simplicial variety (see \cite{D3} 6.1) it is equipped with a
mixed Hodge structure. This structure is not only
pure, but also of Hodge type (i.e. $(k,k)$ in $2k$-th
cohomology).
Every class is represented by an algebraic cycle.\end{re}

\begin{re}\rm The same statement holds for arbitrary smooth $G$-varieties
consisting of finitely many orbits. In fact $H^*_G(X)$ is generated by
algebraic cycles. We will develop this remark in a subsequent
paper \cite{FW}. We want to keep the present paper as elementary as
possible. Therefore we do not treat the case of general $G$-varieties
here. This would demand at least relying on the theory of mixed Hodge
structures (in the smooth case). We prefer to exploit the map $\phi_p$.

\end{re}

\noindent {\it Proof of \ref{pec}.} Our proof is the induction on the
number of cones.

\noindent 1) If $\Delta$
consists of one cone $\sigma$, then
$$H^*_T(X_\Delta)=H^*(B(T_\sigma))\,,$$
where $T_\sigma$ is the subtorus with Lie algebra spanned by
$\sigma$. The isomorphism is induced by the inclusion of the minimal
orbit $T/T_\sigma\mono X_\Delta$. The action of $\phi_p$ on
$H^*(B(T_\sigma))$ is as desired.

\noindent 2) Now suppose that $X$ is decomposed into a
sum of toric varieties $X=X_1\cup X_2$. The Mayer-Vietoris sequence
splits into short exact sequences:
$$0\rightarrow H_T^{2k}(X)\mono H_T^{2k}(X_1)\oplus
H_T^{2k}(X_2)\;\epii\, H_T^{2k}(X_1\cap X_2)\rightarrow 0\,.$$
The groups $H_T^{2k}(X_i)$ are pure, therefore the subgroup of their direct
sum is pure too.\qed

\begin{re}\rm From the intuitive point of view Theorem \ref{pec} is
clear: If one identifies the equivariant cohomology with the
piecewise polynomial functions (Stanley-Reisner ring in the
complete case), then the gradation is given by the homogeneity
degree. Purity means that homogeneous functions are
homogeneous.\end{re}

\subsection{Singular case}
The construction of the equivariant intersection cohomology in
terms  of  the  fan  is  more  complicated, $IH^*_T(-)$ can  be described
axiomatically
(\cite{BBFK1}, Def.~3.1). As in the smooth case the odd part vanishes.
The even part is equipped with an action of $\phi_p$. The map
$IH^*_T(X)\rightarrow IH^*(X)$ is $\phi_p$-invariant. 
We generalize Theorem \ref{pec}:

\begin{th} \label{peic} The equivariant intersection cohomology of a
toric variety is pure.\end{th}

\begin{re}\rm Proceeding as in \cite{BJ} we can prove the same
statement for arbitrary $G$-variety, provided, that it consists of
finitely many orbits and singularities are not to bad.
This is the case for spherical varieties. See \cite{FW}.\end{re}

\Proof The proof of the theorem is analogous except the first step.
We will introduce an induction on $\dim X$. If the dimension is one,
then $X_\Delta$ is smooth ($X_\Delta=\Pt^1$, $\Ct$ or $\Ct^*$) and
the equivariant intersection 
cohomology is pure.

We recall the definition of equivariant formality:

\begin{df} We say that $X$ is equivariantly formal if $IH^*_T(X)$ is a free 
module
over $H^*(BT)$. \end{df}

Equivariant formality is equivalent to the degeneration on $E_2$ of
the spectral sequence:

\begin{eq}\hfil\label{fss} $E^{k,l}_2=H^k(BT)\otimes IH^l(X)\Rightarrow
IH^{k+l}_T(X)\,.$\end{eq}

\noindent In the toric case it is also equivalent to vanishing
of $IH^*(X)$ in odd degrees. Other equivalent conditions are stated in
\cite{BBFK1}, Lemma 4.1 and \cite{BBFK2}, Theorem 3.8.

For the inductive step we will need the following lemma:

\begin{lem} \label{le} Assume that $X$ is equivariantly formal. Then
$IH^*_T(X)$ is pure if and only if $IH^*(X)$ is pure.\end{lem}

\Proof "$\Rightarrow$" If $IH^*_T(X)$ is equivariantly formal then the map
$IH^*_T(X)\rightarrow IH^*(X)$ is surjective. Purity
of $IH^*(X)$ follows. "$\Leftarrow$" If $IH^*(X)$ is equivariantly
formal, then the terms of the spectral sequence \ref{fss} are pure.
Therefore $IH_T^*(X)$ is pure.
\qed
\vskip8pt

\noindent {\it Proof of \ref{peic} cont.} 1) If $\Delta$ consists of a
single cone, then (we divide $X_\Delta$ by a finite group if
necessary, as in the proof Theorem 4.4 in \cite{BBFK1}) there is a
decomposition
$$X_\Delta\simeq T/T_\sigma\times C(X_{L\sigma})\,,$$
where $L_\sigma$ is a fan in ${\bf t}_\sigma/{\rm
lin}(\alpha)$ with $\alpha\in int(\sigma)$. The toric variety
$X_{L\sigma}$ comes with an ample bundle and $C(X_{L\sigma})$ is the
affine cone over $X_{L\sigma}$. It is the toric variety associated
to the cone $\sigma$ considered in ${\bf t}_\sigma$. The variety
$X_{L\sigma}$ is complete, thus it is equivariantly formal (\cite{GKM}). The
equivariant intersection cohomology is pure by the inductive
assumption on the dimension of $X$, so $IH^*(X_{L\sigma})$ is pure by
\ref{le}.
The intersection cohomology of the cone is the primitive part of
$IH^*(X_{L\sigma})$:
$$IH^{n-i}(C(X_{L\sigma}))
=\left\{\matrix{ker(IH^{n-i}(X_{L\sigma})\mpr{h^{i+1}}IH^{n+i+2}(X_{L\sigma}))
&\quad{\rm for}\; i\geq 0\cr
0\hfill& \quad{\rm for}\; i< 0&.\cr}\right.$$
Here $h$ is the class of the hyperplane section. Again we see that
$IH^*(C(X_{L\sigma}))$ is pure. The cone $C(X_{L\sigma})$ is equivariantly 
formal
(since the odd part vanishes). Therefore $IH^*_T(C(X_{L\sigma}))$ is
pure by \ref{le}.

2) The induction with respect to the number of cones is the same as in
the proof of \ref{pec}. One should remember that at this point
decomposition theorem (\cite{BBD}) is used to justify that
$IH^{odd}_T(X)=0$. \qed

\begin{re}\rm Proof that the odd equivariant cohomology vanish and
that even equivariant cohomology is pure may be done in one
shot.\end{re}

\begin{cor} \label{purpur0} Intersection cohomology of $X$ is
pure if and only if $X$ is equivariantly formal.\end{cor}

\Proof The vanishing of $IH^{odd}(X)$ implies degeneration of the
spectral sequence \ref{fss}. To prove the converse
suppose $X$ is equivariantly formal. Equivariant intersection cohomology is 
always
pure by \ref{peic}. Purity of $IH^*(X)$ follows from \ref{le}.\qed

We can restate \ref{purpur0}:

\begin{cor} \label{purpur} Intersection cohomology of $X$ is pure if
and only if it vanishes in odd degrees.\end{cor}

Warning: Our corollaries do not hold for an arbitrary algebraic
variety with a torus action.

\section{Koszul duality}
Let $X$ be a topological space acted by a torus $T$. Borel
construction produces a space and a map:
$$X_T=ET\times_TX\mpr{} BT\,.$$
The homotopy type of $X$ can be recovered from $X_T$
by a pull-back diagram:
$$\matrix{X&\sim&ET\times X&=&ET\times_{BT}X_T&\mpr{}&X_T\cr
&&&&\mpd{}&\p&\mpd{}\cr
&&&&ET&\mpr{}& BT&,\cr}$$
The map $ET\times X\rightarrow X_T$ is a fibration with the
fiber $T$. In the stack language it is the quotient map $X\rightarrow
X/T$, \cite{Si}.

\begin{re}\rm The $T$-homotopy type of $X$ is not preserved by the
procedure 
$$X\mapsto X_T\mapsto ET\times_{BT}X_T=ET\times X$$
since
the action of $T$ on $ET\times X$ is free and the action on
$X$ does not have to be.\end{re}

On the level of cohomology one has a spectral
sequence
\begin{eq} \hfil \label{dss} $E^{k,l}_2=IH^k_T(X)\otimes H^l(T)\Rightarrow
IH^{k+l}(X)\,.$ \end{eq} 
(Note that $X_T$ is simply-connected.) To
recover $IH^*(X)$ one would have to know a complex defining
equivariant intersection cohomology: precisely an element of the
derived category of $H^*(BT)$-modules. This
is a form of Koszul duality as studied in \cite{GKM}, \cite{AP},
\cite{MW}, \cite{Fr1}. 

\begin{re}\rm After certain renumbering of the entries the $E_r$ table of the
spectral sequence \ref{dss} are isomorphic to $E_{r-1}$ of the
Eilenberg--Moore spectral sequence (\cite{EM, Sm}). The analog of
the Eilenberg--Moore spectral sequence for intersection cohomology is
studied in \cite{FW}. It is shown there that the weight
structure is inherited.\end{re}

In the toric case the spectral sequence is acted by $\phi_p$. The
cohomology of $T$ is not pure: $H^l(T)=\Lambda^l{\bf t}^*$ is of
weight $2l$. Because (by \ref{peic}) the equivariant intersection cohomology 
is
pure, the weight of $E^{k,l}_2$ is $k+2l$. The
differentials $$d_{(2)}:E^{k,l}_2\rightarrow E^{k+2,l-1}_2$$ do
preserve weights.

\begin{pr} The higher differentials
$$d_{(i)}:E^{k,l}_i\rightarrow E^{k+i,l-i+1}_2$$
for $i>2$ vanish.\end{pr}

\Proof The weight of the target is $k+i+2(l-i+1)=k+2l-i+2<k+2l$ for
$i>2$. Hence $d_{(i)}$ vanishes.\qed

The differential $d_{(2)}$ has the usual Koszul form:
$$d_{(2)}(x\otimes \xi)=\sum_{j=1}^n x\lambda^j\otimes
i_{\lambda_j}\xi\,, \quad{\rm for}\; x\in IH^*_T(X)\,,\;\;\xi\in\Lambda{\bf
t}^*\,.$$
Here $\{\lambda_j\}$ is the basis of $\bf t$ and
the elements $\lambda^j$ of the dual basis are identified with generators of
$H^*(BT)=S{\bf t}^*$. We obtain a description of the nonequivariant
intersection cohomology:

\begin{th}\label{kosz} Intersection cohomology of a toric variety
is isomorphic to the cohomology of the Koszul complex
$IH^*_T(X)\otimes \Lambda{\bf t}^*$. The filtration induced from the
gradation of $\Lambda{\bf t}^*$ is the weight filtration up to a
shift.\end{th}

\begin{re}\rm The analogous statement holds for arbitrary $G$
varieties with $H^*(G)$ instead of $\Lambda{\bf t}^*$,
provided that $IH_G^*(X)$ is pure.
The proof will appear in \cite{FW}.\end{re}

\begin{re}\rm In the simplicial case the complex $H^*_T(X)\otimes
\Lambda{\bf t}^*$ has appeared in \cite{BP}, 4.2.2 and \cite{Fr1}.
Ishida complex (\cite{Od}) is related to it.
If $X$ is simplicial toric variety, then our complex contains and is
quasiisomorphic to Ishida complex.
\end{re}

Let us describe the weight filtration precisely. The Koszul complex
 splits into a direct sum of
subcomplexes
\begin{eq}\label{Cel}\hfil
$IH^*_T(X)\otimes \Lambda{\bf t}^*=
\bigoplus_{l=0}^{2n}C_{[l]}^*\,,\qquad
C_{[l]}^i=IH_T^{2(k-l)}\otimes\Lambda^{2l-k}{\bf t}^*\,.$\end{eq}
The weight filtration of intersection cohomology also splits into the
sum of eigenspaces of the Frobenius action
$$Gr^W_{2l}IH^k(X)\simeq (IH^k(X))_{p^l}=H^k(C^*_{[l]})\,.$$
Although the complexes $C^*_{[l]}$ might be nonzero in high degrees,
their cohomology vanish:
$$H^k(C^*_{[l]})=0 \quad {\rm for}\; k>l\; {\rm or}\; k>2n\,.$$
In particular we have:
\begin{cor} The pure component of $IH^k(X)$ is equal
to the cokernel
of $d_{(2)}:IH^{k-2}_T(X)\otimes {\bf t}^*\rightarrow IH^k_T(X)$.
Equivalent descriptions of $W_kIH^k(X)$ are
$$IH_T^*(X)/{\bf m}IH_T^*(X)=IH_T^*(X)\otimes_{S{\bf t}^*}\Ct\,,$$
where ${\bf m}= S^{>0}{\bf t}^*$ is the maximal ideal.
\end{cor}

\begin{ex} \rm Let us come back to the Example \ref{prz}. It is the
toric variety associated to the fan consisting of two quarters in
$\Rt^2$ having only the origin in common. According to \S\ref{poly}
the equivariant cohomology can be identified with the pairs of
polynomials in two variables having the same value at the origin.
The spectral sequence \ref{dss} has the $E_2$ table
\vskip6pt
{\begin{tabular}{c|cccccccccc}
2&$\Ct(2)$&0&$\Ct^4(3)$&0&$\Ct^6(4)$&0&$\Ct^8(5)$&0&$\Ct^{10}(6)$&\dots\\
1&$\Ct^2(1)$&0&$\Ct^8(2)$&0&$\Ct^{12}(3)$&0&$\Ct^{16}(4)$&0&$\Ct^{20}(5)$&
\dots\\
0&$\Ct$&0&$\Ct^4(1)$&0&$\Ct^6(2)$&0&$\Ct^8(3)$&0&$\Ct^{10}(4)$&\dots\\ \hline
&0&1&2&3&4&5&6&7&8&\\
\end{tabular}}
\vskip6pt
\noindent The $E_3$ table is
\vskip6pt
{\begin{tabular}{c|cccccccccc}
2&0&0&0&0&0&0&0&0&0&\dots\\
1&0&0&$\Ct(2)$&0&0&0&0&0&0&\dots\\
0&$\Ct$&0&$\Ct^2(1)$&0&0&0&0&0&0&\dots\\ \hline
&0&1&2&3&4&5&6&7&8&\\
\end{tabular}}
\end{ex}

\begin{re}\rm The dual spectral sequence \ref{fss} does
not have to degenerate on $E_3$. The counterexample is
$\Ct^2\setminus\{0\}\sim S^3$ with the action of
$\Ct^*\times\Ct^*$ on the coordinates. The tables $E_2$, $E_3$ and
$E_4$ are equal. The differential $d_{(4)}$ is nontrivial:
$d_{(4)}(1\otimes x)=\lambda^1\lambda^2\otimes 1$, were $x\in
H^3(\Ct^2\setminus\{0\})$ is the generator and $\lambda^1$,
$\lambda^2$ is the basis of ${\bf t}^*$.
\end{re}

\section{Poincar\'e polynomials}

For a toric variety we define the virtual Poincar\'e polynomial by
the formula:
$$IP_{cld}(X)=\sum_{k,l}(-1)^k \dim Gr_{2l}^W IH^k(X)q^l\,.$$
The subscript $cld$ indicates that the supports of cohomology are
closed, whereas the usual Poincar\'e polynomials are taken with
respect to compact supports.
One can find the same formula for cohomology with compact supports
in \cite{Fu} \S4.5 and for intersection cohomology in \cite{DL}
\S5.5. Our approach does not refer to the abstract theory.
The polynomial $IP_{cld}(X)$ can be treated as Euler
characteristic with weights:
$$IP_{cld}(X)=\chi^W(IH^*(X)):=\sum_l\chi(Gr_{2l}^W IH^*(X))q^l\,.$$
Due to purity of equivariant intersection cohomology the equivariant
Poincar\'e polynomial is simpler to define:
$$IP_{cld}(X_T)=\chi^W(IH^*_T(X))=\sum_{l} \dim IH^{2l}_T(X)q^l\,.$$
In general it would be the alternating sum as before.
Poincar\'e polynomials were studied
in \cite{BBFK2}, \S4. The recursive definition of $IP_{cld}(X)$ applies in
our case, as it will be clear from the following propositions. We
note that for complete toric varieties our Poincar\'e polynomial
differs from the one of \cite{BBFK2} by a substitution $q=t^2$,
whereas for noncomplete $X$ the polynomial of \cite{BBFK2} can be
defined as the weighted Euler characteristic of intersection
cohomology with compact supports.

\begin{pr} \label{mult}The Poincar\'e polynomials are related by the formula
$$IP_{cld}(X_T)(1-q)^n=IP_{cld}(X)\,.$$\end{pr}

\Proof The weighted Euler characteristic is multiplicative,
therefore 
$$IP_{cld}(X)=\chi^W(IH^*(X))=\chi^W(IH_T^*(X)\otimes H^*(T))=
\chi^W(IH_T^*(X))\chi^W (H^*(T))\,.$$
The Poincar\'e polynomial of $H^*(T)$ is $(1-q)^n$, thus
$$IP_{cld}(X)=IP_{cld}(X_T)(1-q)^n\,.$$
\qed

Our Poincar\'e polynomials $IP_{cld}(X)$ are not additive with
respect to union of strata, but:

\begin{pr} The polynomials $IP_{cld}(X)$ and $IP_{cld}(X_T)$ are additive 
in the following sense: if $X=X_1\cup X_2$, then
$$IP_{cld}(X)=IP_{cld}(X_1)+IP_{cld}(X_2)-IP_{cld}(X_1\cap X_2)$$
and the same for $IP_{cld}(X_T)$.
\end{pr}

\Proof The additivity of $IP_{cld}(X_T)$ follows from the fact that
Mayer-Vietoris sequence splits into the short exact sequences.
The additivity of $IP_{cld}(X)$ follows from \ref{mult}.
\qed

The recursive formula for $IP_{cld}(X)$ is the following:

$$IP_{cld}(X_\Delta)=\sum_{\sigma\in\Delta}\partial
IP_{cld}(X_{L\sigma})(1-q)^{n-\dim\sigma}\,,$$
where
$$\partial IP_{cld}(X)=\tau_{>{\dim X\over 2}}(1-q)IP_{cld}(X)$$ if
$\dim X\geq 0$ and $\partial IP_{cld}(\emptyset)=1$.
The symbol $\tau_{>{\dim X\over 2}}$ denotes truncating the
coefficients of the monomials of degrees which are not $>{\dim X\over 2}$.
The above formula holds for equivariantly formal $X$ by \cite{BBFK2} and
Poincar\'e duality. Hence it
holds for an affine $X$. By additivity it holds for any $X$.

Andrzej Weber

Instytut Matematyki, Uniwersytet Warszawski,

ul. Banacha 2, 02-097 Warszawa, POLAND

e-mail: aweber@mimuw.edu.pl
\end{document}